\documentclass[10pt,a4paper,twoside]{article}
\usepackage{a4,amsfonts}

\newcommand{\Z}{{\mathbf{Z}}}   
\newcommand{\Q}{{\mathbf{Q}}}
\def\qed{\hspace{\fill}\hbox{${\vcenter{\vbox{              
	\hrule height 0.4pt\hbox{\vrule width 0.4pt height 6pt
	\kern5pt\vrule width 0.4pt}\hrule height 0.4pt}}}$}}

\newtheorem{theorem}{Theorem}
\newtheorem{lemma}{Lemma}[section]
\newtheorem{conjecture}[lemma]{Conjecture}
\newtheorem{definitionhelp}[lemma]{Definition}
\newtheorem{remarkhelp}[lemma]{Remark}
\newtheorem{questionhelp}[lemma]{Question}

\newtheorem{proposition}[lemma]{Proposition}
\newenvironment{definition}{\begin{definitionhelp}\rm}{\end{definitionhelp}}
\newenvironment{remark}{\begin{remarkhelp}\rm}{\end{remarkhelp}}

\def\runninghead#1#2{\pagestyle{myheadings}
\markboth{{\protect\footnotesize\it{\quad #1}}\hfill}
{\hfill{\protect\footnotesize\it{#2\quad}}}}
\input epsf
\newcommand{\epsfig}[1]{\noindent
	\epsfbox{#1}}
\newcommand{\cepsfig}[1]{\par\smallskip\noindent
	\centerline{\epsfig{#1}}\smallskip}
\newcommand{\twocepsfig}[2]{\par\smallskip\noindent
	\centerline{\epsfig{#1}}\par\smallskip
	\centerline{\epsfig{#2}}\smallskip}

\newcommand{\eqnfig}[2]{
	$$\begin{array}{c}
	\epsfig{#1}
	\end{array} 
	\eqno{\rm (#2)}
	$$}
\newcommand{\intextfig}[1]{
	\begin{array}{c}
	\epsfig{#1}
	\end{array} 
	}

\newcommand{\ti}{\hat}
\newcommand{\BB}[2]{{\cal B}_{#1}^{\,#2}}
\newcommand{\B}[2]{\ti {\cal B}_{#1}^{\,#2}}
\newcommand{\F}[1]{{\cal F}_{#1}^{}\,}
\newcommand{\G}[1]{{\cal G}_{#1}^{}\,}
\newcommand{\C}[1]{{\cal C}_{#1}^{}\,}
\newcommand{\epimorph}{{\rightarrow\hspace{-9pt}\rightarrow\;}}

\newcommand{\cc}[1]{\cite{cc#1}}

\begin{document}

\runninghead{\quad Jan A.~Kneissler}{The Ladder Filtration \quad}
\thispagestyle{empty} 
\title{{\large\it On Spaces of connected Graphs III}\\[0.2cm] \smallskip The Ladder Filtration}
\author{Jan A.~Kneissler}
\maketitle
\thispagestyle{empty} 

\vspace{6pt}

\begin{abstract}
{
\noindent
A new filtration of the spaces of tri-/univalent graphs $\BB{m}u$ that 
occur in the theory of finite-type invariants of knots and $3$-manifolds
is introduced.\\
Combining the results of the two preceding articles, 
the quotients of this filtration are modeled by spaces of graphs 
with two types of edges and four types of vertices, and 
an upper bound for $\dim \BB{m}u$ in terms of
the dimensions of the filtration quotients is given.
The degree $m$ up to which $\BB{m}u$ is known is raised by two
for $u=0$ and $u=2$.
}
\end{abstract}

\addtocounter{section}{-1}
\section{Definitions and former results}
This paper is based heavily on results of the preceding articles
\cc1 and \cc2. In this section we briefly repeat the 
relevant definitions and list results of \cc1 that will be used. 

First recall, that in the context of Vassiliev theory 
graphs with vertices of order one and three (uni-, trivalent vertices)
are subjected to some local relations (i.e.~a quotient of a 
vector space with a basis of such graphs is quotiented 
by linear combinations of graphs differing only inside a small region 
in way that is specified by the "relation". Edges that go out off the region
are clipped, but the resulting "ends" should not be confused with univalent 
vertices: all relations in this paper do not involve any univalent vertices.
Furthermore, the ends are supposed to be arranged correspondingly in
all terms of a relations).

The incoming edges at every trivalent vertex are supposed to be 
cyclically ordered and the graphs are always drawn such that these
orderings correspond with the orientation of the plane.
A connected graph with $2m$ vertices will be called {\sl diagram of degree }$m$.

\begin{definition}
\label{defB}
\begin{eqnarray*}
\BB{m}u &:=& \Q\langle \mbox{ diagrams of degree } m \mbox{ with } u
\mbox{ univalent vertices }\rangle\mbox{ } / 
\mbox{ (AS), (IHX)} \\
\B{m}u &:=& \Q\langle \mbox{ diagrams of degree } m \mbox{ with } u
\mbox{ univalent vertices }\rangle\mbox{ } / 
\mbox{ (AS), (IHX), (x)} 
\end{eqnarray*}
where (AS), (IHX) and (x) are the following local relations ((x) is rather 
a family of relations):
\eqnfig{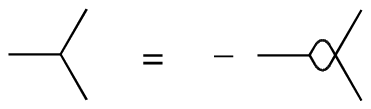}{AS}
\eqnfig{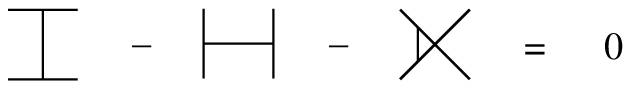}{IHX}
\eqnfig{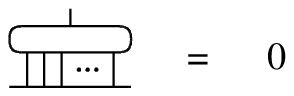}{x}
\end{definition}

\noindent
In \cc1 we have found a set of consequences of the relations in 
$\B{m}u$ that will be repeated here, for the readers convenience:
\eqnfig{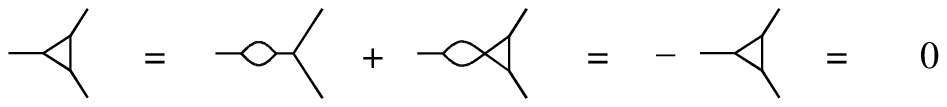}{t}
\eqnfig{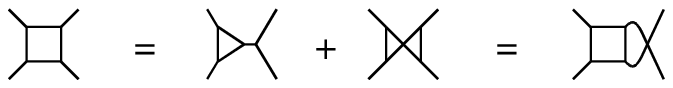}{s}
\eqnfig{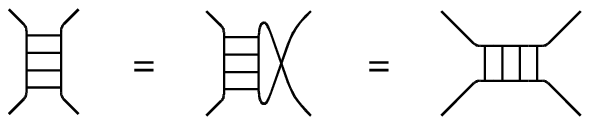}{LS}
\eqnfig{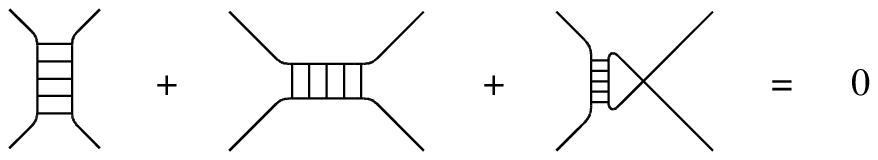}{LIHX}
\eqnfig{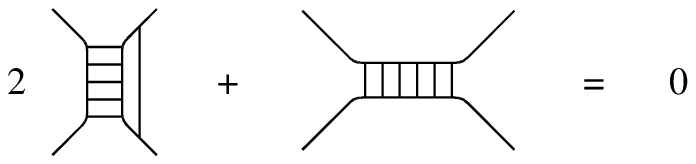}{LI}
\eqnfig{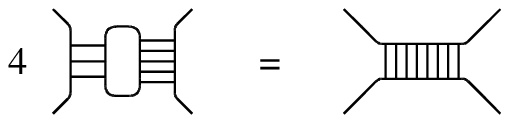}{LL}
The last four relations involve subgraphs of a special form, that will 
be called {\sl ladders}. A ladder with $n\geq 2$ rungs is called $n${\sl -ladder};
it consists of $3n+2$ edges. The four edges connecting the ladder with
the rest of a diagram are called {\sl ends} of the ladder. A 2-ladder 
is also called {\sl square}. An $n$-ladder is called {\sl odd} and {\sl even}
according to the parity of $n$. The set of all ladders of a given diagram 
is partially  ordered by the subgraph relation. 
Ladders that are maximal elements with respect to this relation are called
{\sl maximal} ladders.

Each of the last four equations, involving 
ladders, represents a whole family of relations according 
to the following rules: 
\begin{enumerate}
\item
Each of the example ladders in the pictures above may be replaced by a $n$-ladder
of same parity ($n \geq 2$).
\item
The relation must remain homogeneous (i.e.~each term of a relation must 
get the same total number of rungs).
\end{enumerate}

\noindent
The main result of \cc1 (Theorem 3, "square-tunneling relation") shall also be restated:
\addtocounter{theorem}{-1}
\begin{theorem}
\label{STR}
Moving two rungs from a $p$-ladder to another $q$-ladder inside a diagram 
(thus the new diagram has a $(p-2)$- and a $(q+2)$-ladder) does not change its value in $\B{m}u$ if $p \geq 4$, $q \geq 2$.
\end{theorem}

\section{Introduction}

In the first of this series of articles we 
mainly examined the spaces $\B{m}u$ 
The second article dealt with a sub-algebra $\Lambda_0$ of the algebra $\Lambda$
that has been defined in \cite{Vo}. 
Nevertheless not $\B{m}u$ and $\Lambda_0$  but rather 
$\BB{m}u$ and $\Lambda$ are the spaces that are of significance for 
the theory of finite type invariants of knots and 
of $3$-manifolds, respectively.

So we made life too easy in \cc1 by using the 
additional relation (x)  and 
in \cc2 we possibly missed some part of $\Lambda$ 
since only the elements $x_n$ were considered.

But recall that $\Lambda$ acts on $\BB{m}u$ and that $\Lambda$ is 
isomorphic (as graded vector space) to $\bigcup_{m\geq 1}\BB{m}0$.
So on one side (in \cc1) we ignored the elements of  $\Lambda_0$ 
(of degree $\geq 1$ at least), whereas on the other side 
(in \cc2) we had a somehow ``orthogonal" point of view. 

Thus it is quite natural to look for a synthesis of both approaches that 
``reassembles" $\BB{m}u$ out of the fragments $\B{m}u$ and $\Lambda_0$.

\subsection{The ladder filtration}

The statement of Theorem \ref{STR} can be reformulated as follows:
A square that is part of a $n$-ladder with $n\geq 4$ is called 
{\sl free square} (since it may move around freely across ladders).
The result of removing a free square (i.e.~reducing the number of rungs of 
a ladder by 2) is called a {\sl reduction}. For each diagram one obtains 
after a finite number of reductions a unique diagram, which will be 
called {\sl complete reduction}. Now Theorem \ref{STR} implies that
diagrams of same degree and same complete reduction represent identical 
elements of $\B{m}u$.
In other words, diagrams may be viewed 
as combination of completely reduced diagrams (containing ladders of 
length $2$ or $3$ only) together with a "cloud" of free squares.

We suggest to filter $\B{m}u$ by the number of such free squares
and the number of odd and even ladders, which motivates the 
following Definitions \ref{def1} and \ref{def2}.
\begin{definition}
\label{def1}
For a diagram $D$ and $i\geq2$ let $n_i(D)$ denote the number of maximal ladders of 
length $i$ in $D$. Let $\delta(D)$ be the triple of 
integers counting the number of free squares, odd ladders and even ladders: 
$$\delta(D) \;\;:=\;\; \left(\sum_{i\geq2} 
n_i(D)\cdot\left\lfloor\frac{i-2}2\right\rfloor,\;\; \sum_{i\geq1} n_{2i+1}(D),\;\;
\sum_{i\geq1} n_{2i}(D)\right)$$
\end{definition}

\begin{proposition}
\label{propcond}
Let $D$ denote a diagram of $\B{m}u$ with $\delta(D) = (f,\,o,\,e)$, then 
\begin{eqnarray*}
f,\,o,\,e &\geq& 0\\
o+e = 0 &\Rightarrow& f=0\\
4f+6o+4e &\leq& 2m-u \\
2f+2o+e &\leq& m-u+1.
\end{eqnarray*}
\end{proposition}

\noindent
The first two conditions are obvious, the third states that 
the numbers of vertices in ladders does not exceed the total number
of trivalent vertices in $D$ and the forth is obtained by the 
fact that the cyclomatic number of $D$ is $m-u+1$ and that  
the ladders in $D$ contain $\sum n(i)\cdot(i-1)$ linearly 
independent cycles.

\begin{definition}
\label{def2}
Let $T(m,u)$ denote the subset of $\Z^3$ consisting of all triples $(f,\,o,\,e)$ 
that satisfy the four conditions of Proposition \ref{propcond}.
Let $\prec$ denote the ordering on $T(m,u)$ that is induced by the lexicographical order of $\Z^3$.
Let $\F{f,o,e}\B{m}u$ denote the subspace 
of $\B{m}u$ that is spanned by all diagrams $D$ with $(f,\,o,\,e) \;\preceq\;
\delta(D)$. 
\end{definition}

\noindent
We get a decreasing filtration ($t$ is the maximal element of $T(m,u)$):
$$ \B{m}u \;=\; \F{0,0,0}\B{m}u \;\;\supset\;\; \F{0,0,1}\B{m}u \;\;\supset\;\;
 \ldots \;\;\supset\;\; \F{t}\B{m}u.$$
The aim of this article is to obtain some information about the 
filtration quotients that will be denoted $\G{f,o,e}\B{m}u$.
A strong source for the hope that this filtration 
is helpful for the study of $\B{m}u$ 
(and thus $\BB{m}u$) is the observation that the quotients 
probably are independent of $f$ (up to a shift of indices). We have found evidence for this in many situations, 
which encouraged us to formulate it as a conjecture.
\begin{conjecture}
\label{conjind}
There are isomorphisms $\G{f_1,o,e}\B{m+2f_1}u \;\cong\; \G{f_2,o,e}\B{m+2f_2}u$ for all integers $m,$ $u$, $f_1$, $f_2$, $o$, $e$ 
satisfying $(f_i,\,o,\,e) \in T(m+2f_i,u)$.
\end{conjecture}

\subsection{Feynman graphs}
To represent the elements of $\G{f,o,e}\B{m}u$, we have to specify 
the places where odd and even ladders are situated. 
Therefore, we will introduce graphs that contain vertices of degree $4$ 
and a second type of edges. The similarity to graphs 
used in the physics of elementary particles is so striking 
(but meaningless) that we could not resist calling them Feynman graphs. 
Some authors already name the diagrams of $\BB{m}u$ this way,
(these are also Feynman graphs in our sense, so our terminology extends the old one). 
Nevertheless, to avoid confusion, we will always say ``diagram" when are 
talking about 
elements of $\BB{m}u$, $\B{m}u$ or $\G{f,o,e}\B{m}u$ and 
``Feynman graph" when we are speaking about the spaces $\C{}$ 
that will be defined in this subsection.

\begin{definition}
A {\sl Feynman graph} is a connected multi-graph with two types of edges
called {\sl normal edges} and {\sl photon edges}. The vertices all 
have valency $1, 3$ or $4$. The ends of photon edges must be trivalent
and no more than one photon edge may arrive at a trivalent vertex.
The ends of a photon edge are called {\sl photon vertices} and 
the other trivalent vertices are called {\sl normal vertices}.
At every normal vertex, a cyclic ordering of the three incoming edges
is specified, but at the other vertices, the edges are not ordered in any way.
\end{definition}

\noindent
In pictures, the edges are drawn with straight $\intextfig{}$ 
and wavy $\intextfig{}$ lines. The cyclic ordering at normal edges
is assumed counter-clockwise in all pictures.
The restrictions imply that Feynman graphs may only contain the following four types of vertices:
$$\begin{array}{cccc}
\intextfig{} & \intextfig{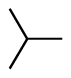} &  \intextfig{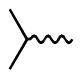} &  \intextfig{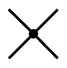} 
\end{array}
$$ 
The {\sl degree} of a Feynman graph is a quadruple of integers 
$(m,\,u,\,p,\,q)$ where $u$ denotes the number of univalent vertices,
$p$ the number of photon edges, $q$ the number of tetravalent vertices and
$2m - u -6p -4q$ is equal to the number of normal vertices ($m$ is always
an integer). Here are two examples of Feynman graphs of
degree $(8,4,1,1)$:
\cepsfig{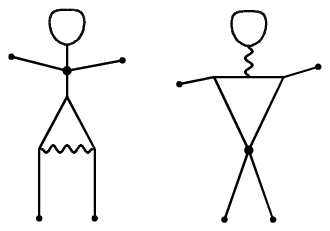}

\noindent
Finally, let $\C{m,u,p,q}$ denote the $\Q$-vector space that is spanned
by all Feynman graphs of degree $(m,\,u,\,p,\,q)$.
The reader might anticipate that the preceding definitions will
ensure the spaces $\G{f,o,e}\B{m}u$ are quotients of these spaces of 
Feynman graphs. In fact, there is an epimorphism
$\varphi_f:\;\C{m,u,o,e} \epimorph \G{f,o,e}\B{m+2f}u$ for all allowed $m,u,f,o,e$
(for simplicity, we denote the all maps with different 
parameters $m,u,o,e$ by the same symbol $\varphi_f$, and hope this will 
not lead to confusion).
\medskip

\noindent
The big question is, how do the kernels of the maps $\varphi_f$ look like?
The best way to answer this question would be to give a complete 
list of local relations that generate the 
kernels for all possible parameters $m,u,f,o,e$. But we do not know 
if such an answer is possible. Furthermore, even if there exists such a list, 
it still might be infinite and almost impossible to find.

At the time, we are only able to present some of the relations we could make out
(see Theorem \ref{theo1}).
They all turn out to be independent of the number $f$,  which strengthens 
our belief in Conjecture \ref{conjind}, which states $\ker\varphi_{f_1}=\ker\varphi_{f_2}$.

\section{Results}

It is not hard to guess how the maps $\varphi_f$ are constructed;
for any Feynman graph $G$ let $r(G)$ denote the diagram that is 
obtained by making the following two substitutions at every tetravalent 
vertex and every photon edge of $G$, respectively:
\begin{eqnarray*}
\intextfig{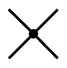} &\rightarrow& \intextfig{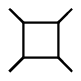} \\
\intextfig{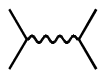} &\rightarrow& \intextfig{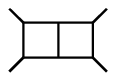}
\end{eqnarray*}
Furthermore, for a diagram $D$ of $\B{m}u$ containing at least one ladder,
let $s(D)$ denote the diagram that is obtained by adding a square to 
an arbitrarily chosen ladder of $D$:
\begin{eqnarray*}
\intextfig{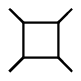} &\rightarrow& \intextfig{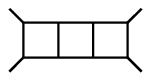}
\end{eqnarray*}
If there are two or more ladders in $D$ then $s(D)$ is not uniquely defined 
by this, but by the square-tunneling relation 
(Theorem \ref{STR}) ensures that different choices yield the same 
element in $\B{m+2}u$.

\begin{theorem}
\label{theo1}
The map $s^f\circ r$ from Feynman graphs to diagrams
extends to a well-defined epimorphism 
$\varphi_f: \; \C{m,u,o,e} \;\epimorph\; \G{f,o,e}\B{m+2f}u$ for all $(f,o,e) \in T(m+2f,u)$.
Elements of the subsequent forms are annihilated by these epimorphisms:
\begin{itemize}
\item
linear combinations that look locally like one of the following:
$$
\intextfig{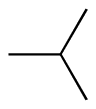} + \intextfig{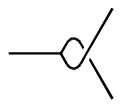}\hspace{-8pt}, 
\eqno{\rm(AS)}
$$
$$
\intextfig{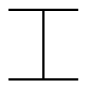} - \intextfig{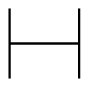} - \intextfig{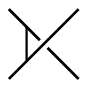}\hspace{-8pt}, 
\eqno{\rm(IHX)}
$$
$$
\intextfig{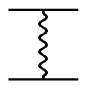} + \intextfig{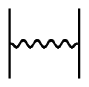} + \intextfig{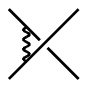}\hspace{-8pt},
\eqno{\rm(LIHX)}
$$
$$
\intextfig{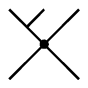} + \intextfig{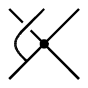} - \intextfig{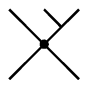} - \intextfig{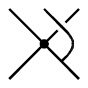}\hspace{-8pt},
\eqno{\rm(4T_a)}
$$
$$
\intextfig{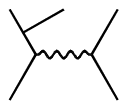} + \intextfig{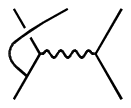} - \intextfig{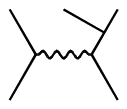} - \intextfig{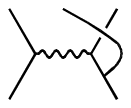}\hspace{-8pt},
\eqno{\rm(4T_b)}
$$
\item
graphs that contain a cycle of length 1,2 (if $m > 2$) or 3,
\item
graphs that contain a cycle of length 4 of one of the following types:
$$
\intextfig{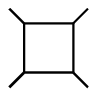},
\intextfig{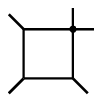},
\intextfig{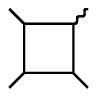},
\intextfig{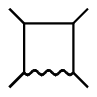},
\intextfig{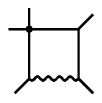},
$$
$$
\intextfig{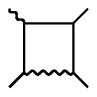},
\intextfig{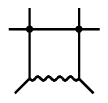},
\intextfig{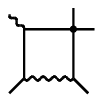},
\intextfig{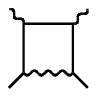},
\intextfig{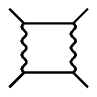},
$$
\item
graphs that contain a cycle of length 5 of one of the following types: 
$$
\intextfig{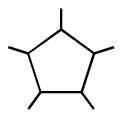},
\intextfig{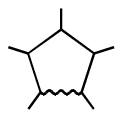},
\intextfig{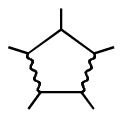}.
$$
\end{itemize}
\end{theorem}

\noindent
The fact that there are a lot of forbidden cycles lets one expect that
$-$ at least in low degrees $-$ the spaces $\G{f,o,e}\B{m}u$ have low
dimensions. For instance, using the Euler characteristics of compact orientable surfaces, it is not hard 
to show that a Feynman graph with genus $g$, $q$ tetravalent vertices and $u < q+6(1-g)$ 
univalent vertices contains a cycle of length $\leq 5$. This implies for $u < 6$ that $\G{0,0,0}\B{m}u$
is spanned by non-planar graphs. In \cite{CD} it has been conjectured that
$\BB{m}0$ is spanned by planar diagrams; if this was true, then 
$\G{0,0,0}\B{m}0 = 0$ for all $m$.

\subsection{Upper bounds for $\dim \BB{m}u$}
We now focus on the spaces $\BB{m}u$ rather than $\B{m}u$.
It is clear that information about the filtration quotients 
$\G{f,o,e}\B{m}u$ will somehow provide upper bounds for $\dim \BB{m}u$. 
With the next statement, we will make this more explicit.

\newcommand{\sqnum}[1]{q_{#1}}
\begin{theorem}
\label{theo2}
Assume that we are given upper bounds for the dimensions of $\G{f,o,e}\B{m}u$
that are independent of $f$, i.e.~there exist numbers $\mu_{m,u,o,e}$ 
satisfying $\dim \G{f,o,e}\B{m+2f}u \;\;\leq\;\; \mu_{m,u,o,e}$ for all $f$.
With the abbreviations $\sqnum{n} := \left\lfloor\frac{n^2}{12}\!+\!\frac12\right\rfloor$
and 
$\delta_{i,j} := \left\{
{1 \;\;\mbox{\scriptsize for }\;i=j} \atop {0 \;\;\mbox{\scriptsize for }\;i\neq j}\right.$,
the following inequality then holds for all $m \geq u \geq 0$:
$$
\dim\BB{m}u \;\;\leq\;\; \sum_{j\,=\,u}^m\left((1+\sqnum{m-j})
\,\mu_{j,u,0,0} \;+\; 
\sqnum{m-j+3} 
\sum_{o+e\geq 1}
\mu_{j,u,o,e}
\right)
\;+\;\sqnum{m-6}\cdot\delta_{u,4}.
$$
The inner sum runs over all $o,e \geq 0$ with $o+e \geq1$, $6o+4e \leq 2m-u$
and $2o+e \leq m-u+1$.
\end{theorem}

\noindent
The idea, of course, is to obtain such upper bounds
by use of Theorem \ref{theo1}, by simply taking
$$\mu_{m,u,o,e} \;\; := \;\;
\dim \big(\; \C{m,u,o,e} \;\;/\;\; \mbox{relations listed in Theorem \ref{theo1}}\;\big).$$

\subsection{Experimental data}

The spaces $\G{f,o,e}\B{m}u$ have been calculated for $m \leq 12$ and $2 \leq u \leq 6$.
This has been done by establishing upper and lower bounds for 
$\dim \G{f,o,e}\B{m}u$. The computation of upper bounds is based on 
Theorem \ref{theo1}; lower bounds are derived by comparing the
statement of Theorem \ref{theo2} with the known exact values of 
$\dim \BB{m}u$ of \cite{Kn}.

As predicted by Conjecture \ref{conjind}, the calculated filtration spaces 
do not depend on $f$ (up to the shift in the degree $m$). All but $15$ of the spaces
in question are trivial. In the non-trivial cases, $\G{f,o,e}\B{m}u$
has dimension one. In the following table we list these spaces together with 
a Feynman graph that is mapped by $\varphi_f$ to a generator (the case $u=0$ which 
can be derived from $u=2$ is also shown).
The non-triviality of these diagrams (as elements of 
$\G{f,o,e}\B{m}u$) has been verified using 
weight systems that come from the Lie-algebras $\mathfrak gl$ and $\mathfrak so$. 

$$\begin{array}{c|c|c|c|c|c|c}
\G{0,0,0}\B10 & 
\G{0,0,0}\B22 & 
\G{f,0,1}\B{4+2f}4&
\G{0,0,0}\B66 & 
\G{f,0,2}\B{7+2f}6&
\G{f,1,0}\B{8+2f}6&
\G{f,2,0}\B{10+2f}6\\
&& \scriptstyle{(0\leq f\leq4)} && \scriptstyle{(0\leq f\leq2)} & \scriptstyle{(0\leq f\leq2)} & \scriptstyle{(0\leq f\leq1)}\\
\hline 
& & & & & &  \\[-7pt]
\hspace{-0.2cm}\intextfig{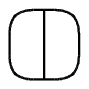}\hspace{-0.5cm} &
\hspace{-0.2cm}\intextfig{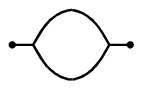}\hspace{-0.5cm} &
\hspace{-0.2cm}\intextfig{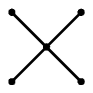}\hspace{-0.5cm} &
\hspace{-0.2cm}\intextfig{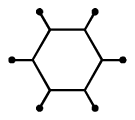}\hspace{-0.5cm} &
\hspace{-0.2cm}\intextfig{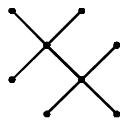}\hspace{-0.5cm} &
\hspace{-0.2cm}\intextfig{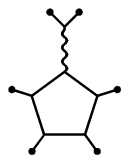}\hspace{-0.5cm} &
\hspace{-0.2cm}\intextfig{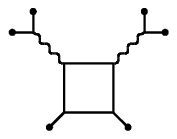}\hspace{-0.5cm} 
\end{array}
$$


\noindent
\begin{remark} 
The upper bounds for $\dim \BB{m}u$ of Theorem \ref{theo2} 
are sharp for $m \leq 12$.
\end{remark} 
Considering the results of \cite{Kn}, one might conjecture that 
$\B{m}u$ is spanned by diagrams the following type (i.e.~caterpillars with
an even number of legs on each of its body segments):
\cepsfig{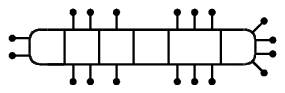}

\noindent
This conjecture is quite daring and we have the feeling that it will turn out 
false. Nevertheless let us mention that it implies the following weaker 
statement, that is more likely to be true:

\begin{conjecture}
\label{conjtriv}
For $m > u+1 > 0$ and $m > 5\left(\frac{u}2-1\right)$
the filtration quotients $\G{f,o,e}\B{m+2f}u$ are trivial for all possible 
parameters $f,o,e$.
\end{conjecture}

\noindent
If Conjecture \ref{conjtriv}
is true then, due to Theorem \ref{theo2}, $\dim \BB{m}u$ grows at most quadratically 
with $m$ for fixed $u$.
To be more specific: together with Theorem \ref{theo2} and 
taking into account the data of the table above, 
Conjecture \ref{conjtriv} implies 
\begin{conjecture}
\label{conjineq}
\begin{eqnarray*}
\dim \BB{m-1}0 \;\;=\;\; \dim \BB{m}{2} &\leq& \left\lfloor \frac{(m-2)^2}{12}+\frac32\right\rfloor\\
\dim \BB{m}4 &\leq& \left\lfloor \frac{(2m-7)^2}{24}+1\right\rfloor\\
\dim \BB{m}6 &\leq& \left\lfloor \frac{(2m-11)^2}{12}+\frac32\right\rfloor.
\end{eqnarray*}
\end{conjecture}
These conjectural inequalities are fulfilled and even sharp 
for $m \leq 12$. 
\medskip

\begin{remark}
The restriction of Conjecture \ref{conjtriv} to the case $u=0$ 
(or equivalently $u=2$) is equivalent to Conjecture 1.2 of \cc2
(claiming $\Lambda=\Lambda_0$, i.e.~Vogel's algebra $\Lambda$ is 
generated by the elements $x_n$ with $n \geq 1$ already, see 
also \cite{Vo}). It is known that $\Lambda=\Lambda_0$ for degrees
$\leq 12$ due to the following
\end{remark}

\noindent
{\bf Computational result:} $\quad$ $\dim\,\G{f,o,e}\B{m}2 = 0\;$
for $\;3 \leq m \leq 14$.
\medskip

\noindent
This also implies that the first inequality of Conjecture \ref{conjineq}
is correct and sharp up to $m = 14$.
The algorithm of the computer program that checked triviality 
for the previously unknown degrees $m=13$ and $m=14$ is based on
the relations of Theorem \ref{theo1} and the fact that in all considered
degrees it suffices to take into account only Feynman graphs of a 
very special form (graphs in which there is a path going 
through all non-normal vertices and through all photon edges, 
not containing any normal vertices). 

\begin{remark}
In chapter 7 of \cite{diss}, using results of Petersen and Chv\'atal,
$\dim\,\G{f,o,e}\B{13}2 = 0\;$ is shown "by hand", i.e.~without aid of a
computer program. Petersen's theorem states that
any $(2n)$-regular graph has a $2$-factor; Chv\'atal's theorem 
gives a sufficient condition for a graph being path-hamiltonian 
in terms of its degree sequence.
\end{remark}

\section{Proof of Theorem \ref{theo1}}

The map $r$ replaces tetravalent vertices by squares and photon edges 
by ladders of length $3$. Because of relation (s), the result does not depend 
(in $\B{m}u$) on the different possible choices of permutations of the four 
incoming edges. 

Due to the square-tunneling relation (Theorem \ref{STR}), adding squares 
to existing ladders is a well-defined operation in $\B{m}u$. Furthermore, the degree of Feynman 
graphs has been chosen to ensure that $r$ is degree-preserving with
respect to the $m$- and $u$-components of the degree and $s$ increases 
the total degree $m$ by $2$. 
So $s^f\circ r$ extends to a homomorphism from $\C{m,u,o,e}$ 
to $\G{f,o,e}\B{m+2f}u$
(there is no problem with the fact that $s$ is only defined
for diagrams with ladders, because $o=e=0$ implies $f=0$, i.e.~$s$ is never 
applied on diagrams without ladders).

Obviously, every diagram of $\G{f,o,e}\B{m+2f}u$ is the image of 
some Feynman graph in $\C{m,u,o,e}$ (to find the Feynman graph, 
simply take the complete reduction of the diagram and invert $r$).

\medskip
\noindent
So the only thing that remains to prove Theorem \ref{theo1} is to
verify the list of relations.
The relations (AS), (IHX), (LIHX) are mapped by $\varphi_f$
to the relation in $\B{m+2f}u$ with the same names; 
(4T$_{\rm a}$) and (4T$_{\rm b}$) are special cases of the fact 
that, taking the (consistently signed) sum of all possible ways of 
attaching an edge at all entries of a region without univalent 
vertices, one obtains a trivial element of $\BB{m}u$, which is 
obviously also trivial in $\B{m}u$ (see Lemma 3.1 a) of \cc1 for 
more details).

\subsection{Loops and double edges}
There are three possible situations with cycles of length $1$:
$$
\intextfig{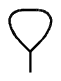},
\intextfig{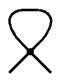},
\intextfig{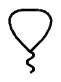}.
$$
The first is trivial by (AS), the last two are trivial because
of the relation (x).
Next, we consider cycles of length $2$ (i.e. double edges):
$$
\intextfig{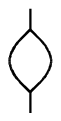},
\intextfig{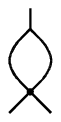},
\intextfig{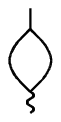},
\intextfig{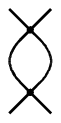},
\intextfig{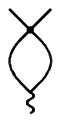},
\intextfig{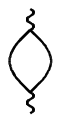},
\intextfig{}.
$$
In the first case, if there is at least one additional trivalent vertex in the 
diagram, we may apply relation (x). Only in the cases 
$m=1,u=0$ and $m=2,u=2$, there exist diagrams having a cycle of the first 
type, but no other trivalent vertices (see the first two columns of the table on the preceding 
page). This explains why we have to make
the requirement $m > 2$. 

The cases 2, 3 are killed by relation (t) and cases 4, 5, 6 are mapped by $\varphi_f$
to diagrams with more than $f$ free squares, because the length of two
ladders are added. In the last case we may apply relation (x).

\subsection{The flyping-trick}

In many situations we apply a useful procedure that we call {\sl flyping}.
It can generally be described like this: take the cycle out of the plane,
make a half twist around an appropriate axis and put it back into the plane. 
By (AS) each normal vertex causes a change of sign. 
Then use relations (AS), (IHX), (4T$_{\rm a}$), (4T$_{\rm b}$) and (LIHX) to take
the cycle into its original form, modulo simpler cycles. 
If the result has opposite sign, then the described operation yields an equation that allows
to express the cycle in terms of simpler diagrams.

We illustrate this trick for two examples of triangles that in this way can be
expressed by diagrams with double edges (relations (LIHX) and (4T$_{\rm a}$) are used):
\twocepsfig{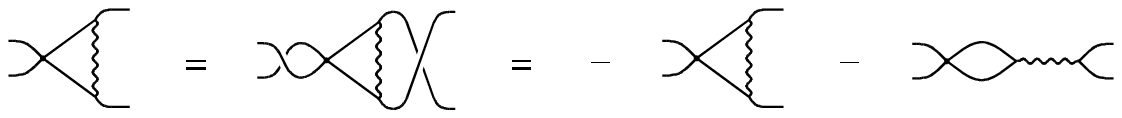}{}

\subsection{Triangles}
\label{sectboxes}
There are thirteen possible cycles of length $3$; let us first discuss the easy ones:
\begin{eqnarray*}
&
\intextfig{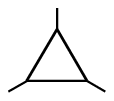},
\intextfig{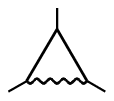},
\intextfig{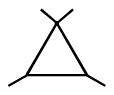},
\intextfig{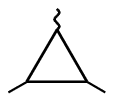},
\intextfig{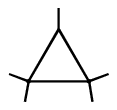},
&\\&
\intextfig{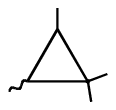},
\intextfig{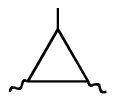},
\intextfig{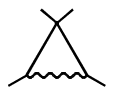},
\intextfig{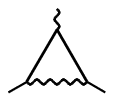}.
&\end{eqnarray*}
The first is killed by relation (t), the second by (x). In the third and forth situations,
the length of the ladder may be increased by including the edge between the two normal vertices.
In the last five cases one may successfully apply the flyping-trick.
Now we have to take a look at the tricky triangles:
$$
\intextfig{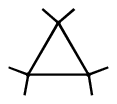},
\intextfig{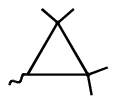},
\intextfig{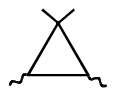},
\intextfig{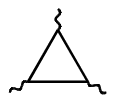}.
$$

\noindent
Let us assume that all free squares are inserted at the ladder on the left side,
when $\varphi_f$ is applied to these graphs.
Then the ladders in the upper and right corner will have length $2$ or $3$.
Here we need a result of \cc1 (equation (37) in section 4.4), which states 
the following local relation:
\eqnfig{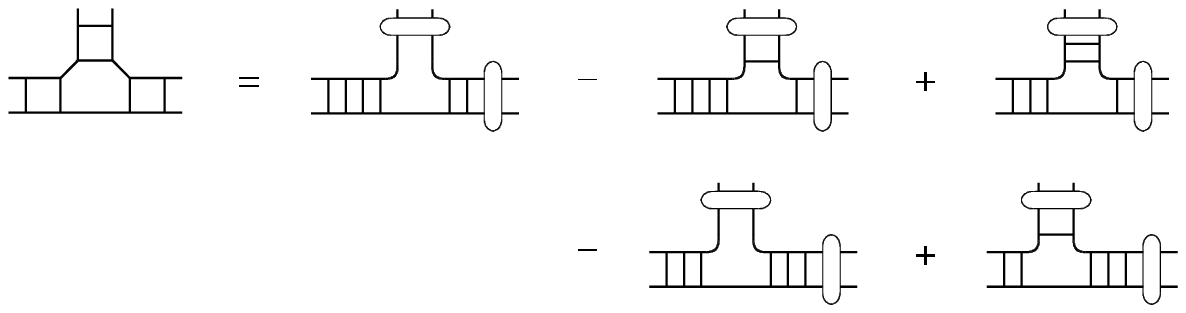}{\star}
The diagrams with rounded boxes shall be interpreted as linear combination of 
4 diagrams by expanding each box the following manner:
\cepsfig{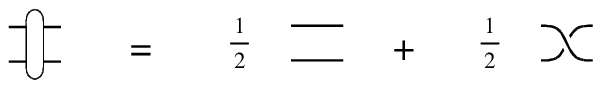}

\noindent
We may apply ($\star$) in all four cases and obtain diagrams in which the
triple (number free squares, number of odd ladders, number of even ladders)
is lexicographically greater than in the original diagram.

\subsection{Simple cycles of length $4$ and $5$}
Again we start with the simple cases:
$$
\intextfig{fig41.ps},
\intextfig{fig42.ps},
\intextfig{fig43.ps},
\intextfig{fig412.ps},
\intextfig{fig411.ps}.
$$
The first contains a square that is not coming from a tetravalent vertex,
the next three can be handled with the flyping-trick and the last is due
to the relation (LL).

%
There are two $5$-cycles that can be reduced by the flyping-trick to the first
and last of the $4$-cycles above:
$$
\intextfig{fig51.ps},
\intextfig{fig55.ps}.
$$

\subsection{Pentagons with a photon edge}
Perhaps the most surprising (and hardest to prove) statement is that 
cycles of this type (i.e.~containing a photon edge and three normal vertices)
are annihilated by $\varphi_f$:
$$
\intextfig{fig52.ps}.
$$

\noindent
To show this, we start with the observation that the occurrence of one of the 
following four situations allows to produce a ladder of length $4$:
$$
\mbox{a)} \intextfig{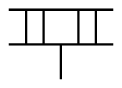}, \quad \quad
\mbox{b)} \intextfig{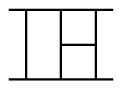}, \quad \quad
\mbox{c)} \intextfig{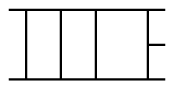}, \quad \quad
\mbox{d)} \intextfig{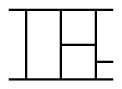}.
$$
In case a) and c) this is done by the flyping trick. Case b) corresponds to
relation (LI) and d) can be reduced to a trivial diagram (because of relation (t))
and type c) and b) in the following way:
\cepsfig{}

\noindent
Let us recall a result of \cc1 (equation (38), section 4.4, it involves a
parameter $n\geq 1$), which 
looks for $n=1$ like this:
\cepsfig{}
If we glue the upper left ends together like this
$\epsfxsize=4cm\intextfig{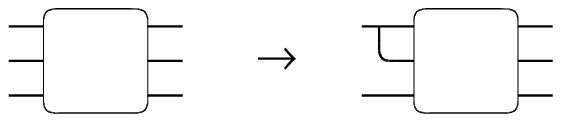},$
we obtain 
$$\intextfig{}.$$
By the sign $\equiv$, we mean ``modulo terms containing a $4$-ladder".
This congruence is true, because three of the diagrams in the middle are of
type a), b) and d).
On the other hand, if we use (IHX) and neglect another diagram of type 
d), we get
$$\intextfig{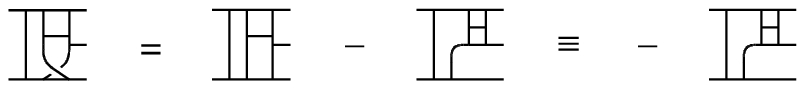}.$$
One may repeat this argument for arbitrary odd ladders instead of the $3$-ladder
(using the general version of equation (38) of \cc1). Thus we
have shown that the following two Feynman graphs are mapped to equivalent 
diagrams in the corresponding filtration quotient:
\cepsfig{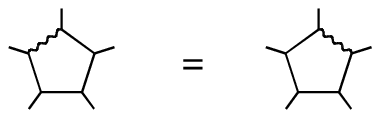}
This finally leads to the following (oval boxes to be interpreted as in section \ref{sectboxes}):
$$\intextfig{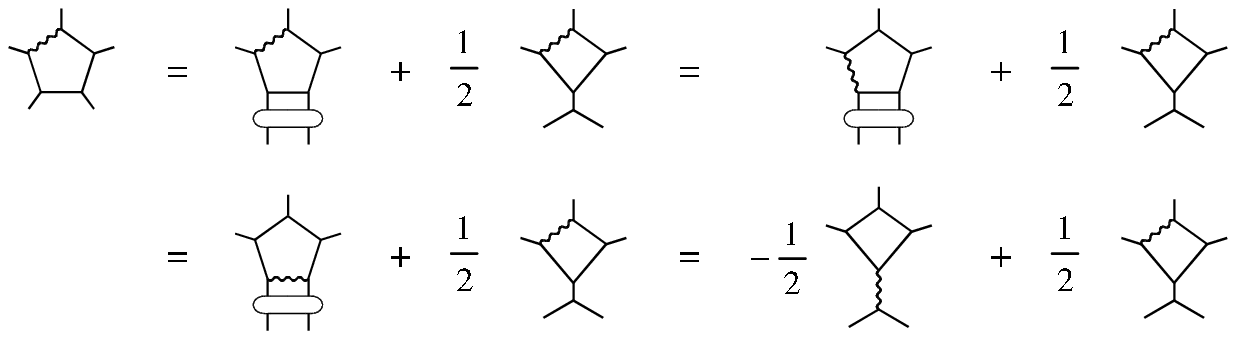}.$$
As we have seen before, these two $4$-cycles lie in the kernel of $\varphi_f$.

\subsection{Cycles with a photon edge and non-normal vertices}
In the preceding subsection we have seen that the occurrence of a pentagon with
a photon edge allows to increase the number of free squares.
This immediately applies in the following 
two situations (we have to sacrifice ladders when reducing the 
pentagon, but this is acceptable since
we gain free squares):
$$
\intextfig{fig49.ps},
\intextfig{fig410.ps}.
$$

\noindent
The only three remaining cycles in our list are the following:
$$
\intextfig{fig46.ps},
\intextfig{fig47.ps},
\intextfig{fig48.ps}.
$$
They all can be reduced by the flyping-trick, if we realize that 
\cepsfig{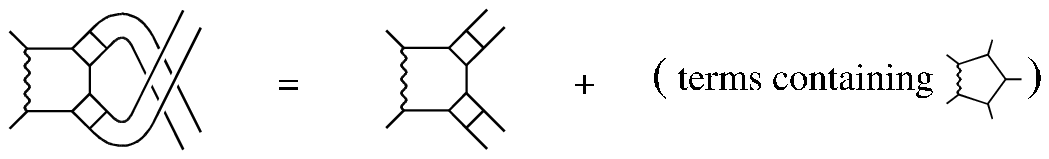}
\noindent
and similar equations with one or both of the squares replaced by $3$-ladders.
\qed

\section{Proof of Theorem \ref{theo2}}


Recall that $\BB{m}u$ may be viewed as $\Lambda_0$-module 
(see \cc2 for a definition of $\Lambda_0$ and its elements $x_n$).
$\B{m}u$ is obtained by quotienting $\BB{m}u$ by all elements that are multiples
of $\Lambda_0$-elements of degree $\geq 1$.
So by Theorem \ref{theo1}, $\BB{m}u$ is spanned the following set: 
($\G{n}\Lambda_0$ denotes the degree $n$ part of $\Lambda_0$)
$$\{\;\lambda\;\varphi_f(d) \;\;\vert\;\; \lambda \in \G{i}\Lambda_0,\;\; d \in \C{m-2f-i,u,p,q},\;\;
(f,p,q)\in T(m-i,u)\;\}.$$ 
In other words: we obtain a spanning set for $\BB{m}u$ by taking  
basis elements of $\C{m^\prime,u,p,q}$ with $m^\prime \leq m$
and blowing up their degree to the correct value by inserting free squares
and elements of $\Lambda_0$ in all possible ways.

To get the desired inequality, we still have to reduce the number of 
generators a little bit. We therefor have to use the fact
that $\BB{m}u$ is not a free $\Lambda_0$-module, what has been 
discovered in \cite{Kn}. 
\medskip

\noindent
By Theorem \ref{STR} we know that the following local
relation holds in $\B{m}u$:
\cepsfig{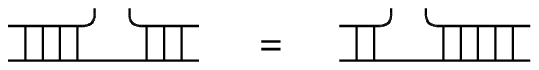}
But in view of Corollary 2.4 of \cite{cc1}, one can hope for 
more: there must also exist a similar local relation in  $\BB{m}u$,
involving elements of $\Lambda_0$ (a local relation
is, in the terminology of \cite{cc1}, nothing else than a trivial
linear combination in $F(6)$; Corollary 2.4 of \cite{cc1} gives
such a trivial element in $\ti F(6)$, the $\Lambda_0$-quotient of $F(6)$).
The corresponding identity in $F(6)$ in fact looks like this:
\cepsfig{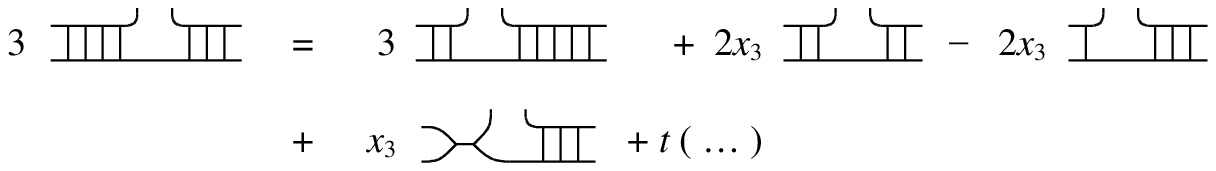}

\noindent
Closing this relation in the following way $\epsfxsize=4cm\intextfig{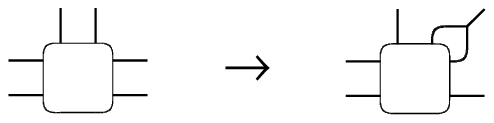}$\hspace{-0.3cm},
we obtain for $n \geq 3$:
$$
x_n\hspace{-0.1cm}\intextfig{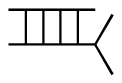}\hspace{-0.2cm} \;\;=\;\; (x_{n+2} + \frac23x_3x_{n-1})\hspace{-0.1cm}\intextfig{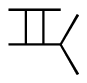}\hspace{-0.2cm} 
-\;\;\frac23x_3x_n\hspace{-0.1cm}\intextfig{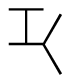}\hspace{-0.2cm} +\;\;\frac13x_3x_n \hspace{-0.1cm}\intextfig{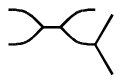}\hspace{-0.2cm} +\;\; t\;(\;\ldots\;)
$$
This means that an element of the form 
$t^ix_n\varphi_f(d)$ with $n \geq 3$ and $f > 0$ can be expressed in terms of
elements containing more $t$-s or a $\Lambda_0$-factor of higher degree.
There is one exception: if $d=\hspace{-0.2cm}\intextfig{}\hspace{-0.45cm}\in\C{4,4,0,1}$ and $f=1$, then there is no
trivalent vertex in $\varphi_1(d)$ other than the eight vertices of the ladder.
We thus may not apply the upper relation to this case.

In all other cases we may lower the number of free squares successively,
as long as there is an $x_n$ with $n\geq 3$ present.
Consequently, $\BB{m}u$ is spanned by elements of the three types
(the third only applies to $u=4$):
\begin{itemize}
\item[a)]
$\lambda\varphi_0(d)$ with $\lambda \in \G{i}\Lambda_0$ and $d \in D_{m-i,u,p,q}$, $(0,p,q)\in T(m-i,u)$,   
\item[b)]
$t^i\varphi_f(d)$  with $f \geq 1$ and $d \in D_{m-i-2f,u,p,q}$, $(f,p,q)\in T(m-i,u)$,
\item[c)]
$\lambda\varphi_1(\hspace{-0.2cm}\intextfig{figcros.ps}\hspace{-0.45cm})$ with $\lambda \in \G{m-6}\Lambda_0$. 
\end{itemize}
To make this argument more convincing, we would have to define 
a filtration of $\BB{m}u$ according to the lexicographical ordering
of the pairs (number of $t$-s, sum of degrees of $x_n$-s with $n \geq 3$), but 
we omit this technical complication.

By Corollary 1.6 of \cc2 we know that $\dim \G{n}\Lambda_0 \leq 1+\sqnum{n}$.
So the number of linearly independent elements of type a) is at most 
$$\sum_{i = m-u}^m\;\;\sum_{p,q}\;(1+\sqnum{i})\dim D_{m-i,u,p,q}.$$

\noindent
The outer sum starts with $i=m-u$, because $D_{m,u,p,q}$ is always trivial 
if $m < u$. The inner sum runs over all allowed $p,q$ i.e.~all
$p,q$ that satisfy $(0,p,q)\in T(m-i,u)$.

\noindent
The number of independent elements of type b) is
$$\sum_{j=u}^{m} \;\;\sum_{p+q\geq1}\;\left\lfloor\frac{m-j}2\right\rfloor\dim D_{j,u,p,q}$$
because $\lfloor\frac{m-j}2\rfloor$ 
is the number of pairs $(i,f)$ satisfying $i\geq 0, f \geq 1, j+i+2f=m$.
The inner sum runs over all allowed $p,q$ with $p + q \geq 1$
(we can only add free squares if there is at least one ladder).

There are $1+q_{m-6}$ elements of type c), but $\lambda=t^{m-6}$ has already
been counted in b). So the additional contribution to the dimension
by type c) elements is at most $q_{m-j}\;\delta_{u,4}$.
\medskip

\noindent
Substituting $j = m-i$ and using $1+\sqnum{n}+\lfloor\frac{n}2\rfloor =
\sqnum{n+3}$, all three contributions sum up to the inequality given in
Theorem \ref{theo2}.
\qed

\begin{verbatim}
 
e-mail:jan@kneissler.info
http://www.kneissler.info
\end{verbatim}

\end{document}